\documentclass [12pt]{article}
\usepackage{amssymb}
\def\lanbox{\hbox{$\, \vrule height 0.25cm width 0.25cm depth 0.01cm \,$}}

\usepackage{color}

\begin{document}

\centerline{\Large\bf On the solvability of
some systems of integro-}

\centerline{\Large\bf differential equations with and without a drift}

\bigskip

\centerline{Messoud Efendiev$^{1,2}$, Vitali Vougalter$^{3 \ *}$}

\bigskip

\centerline{$^1$ Helmholtz Zentrum M\"unchen, Institut f\"ur Computational
Biology, Ingolst\"adter Landstrasse 1}

\centerline{Neuherberg, 85764, Germany}

\centerline{e-mail: messoud.efendiyev@helmholtz-muenchen.de}

\centerline{$^2$ Department of Mathematics, Marmara University, Istanbul,
T\"urkiye}

\centerline{e-mail: m.efendiyev@marmara.edu.tr}

\bigskip

\centerline{$^{3 \ *}$ Department of Mathematics, University of Toronto,
Toronto, Ontario, M5S 2E4, Canada}

\centerline{e-mail: vitali@math.toronto.edu}

\bigskip
\bigskip
\bigskip

\noindent {\bf Abstract.}
We prove the existence of solutions for
some integro-differential systems containing equations with and without the
drift terms in the $H^{2}$ spaces by virtue of the fixed point
technique when the elliptic equations contain second order differential
operators with and without the Fredholm property, on the whole real line or
on a finite interval with periodic boundary conditions.
Let us emphasize that the study of the
system case is more complicated than of the scalar situation and requires to
overcome more cumbersome technicalities.

\bigskip
\bigskip

\noindent {\bf Keywords:} solvability conditions, non-Fredholm
operators, integro-differential systems, drift terms

\noindent {\bf AMS subject classification:} 35J61,\ 35R09, \ 35K57

\bigskip
\bigskip
\bigskip
\bigskip

\section{Introduction}

 \noindent
 We recall that a linear operator $L$ acting from a Banach
 space $E$ into another Banach space $F$ satisfies the Fredholm
 property if its image is closed, the dimension of its kernel and
 the codimension of its image are finite. Consequently, the
 problem $Lu=f$ is solvable if and only if $\phi_i(f)=0$ for a
 finite number of functionals $\phi_i$ from the dual space $F^*$.
 These properties of the Fredholm operators are broadly used in many
 methods of the linear and nonlinear analysis.

 \noindent
 Elliptic equations in bounded domains with a sufficiently smooth
 boundary satisfy the Fredholm property if the ellipticity
 condition, proper ellipticity and Lopatinskii conditions are
 fulfilled (see e.g. \cite{Ag}, \cite{E09}, \cite{LM}, \cite{Volevich}).
 This is the main result of the theory of the linear elliptic problems.
 In the case of the unbounded domains, these conditions may
 not be sufficient and the Fredholm property may not be satisfied.
 For example, the Laplace operator, $Lu = \Delta u$, in $\mathbb R^d$
 fails to satisfy the Fredholm property when considered in
 H\"older spaces, $L : C^{2+\alpha}(\mathbb R^d) \to C^{\alpha}(\mathbb
 R^d)$, or in Sobolev spaces,  $L : H^2(\mathbb R^d) \to L^2(\mathbb
 R^d)$.

 \noindent
 Linear elliptic problems in unbounded domains satisfy the
 Fredholm property if and only if, in addition to the assumptions
 stated above, the limiting operators are invertible (see \cite{V11}). In some
 simple cases, the limiting operators can be constructed explicitly. For
 example, if
 $$
 L u = a(x) u'' + b(x) u' + c(x) u , \;\;\; x \in \mathbb R ,
 $$
 where the coefficients of the operator have limits at the infinities,

 $$ a_\pm =\lim_{x \to \pm \infty} a(x) , \;\;\;
 b_\pm =\lim_{x \to \pm \infty} b(x) , \;\;\;
 c_\pm =\lim_{x \to \pm \infty} c(x) , $$
 the limiting operators are given by:

 $$ L_{\pm}u = a_\pm u'' + b_\pm u' + c_\pm u . $$
 Since the coefficients here are the constants, the essential spectrum of
 the operator, that is the set of the complex numbers $\lambda$ for
 which the operator $L-\lambda$ does not satisfy the Fredholm
 property, can be found explicitly by means of the Fourier
 transform:
 $$
 \lambda_{\pm}(\xi) = -a_\pm \xi^2 + b_\pm i\xi + c_\pm , \;\;\;
 \xi \in \mathbb R .
 $$
 The invertibility of the limiting operators is equivalent to the
 condition that the essential spectrum does not contain the origin.

 \noindent
 In the cases of the general elliptic equations, the same assertions hold true.
 The Fredholm property is satisfied if the essential spectrum does not contain 
 the origin or if the limiting operators are invertible. However,
 these conditions may not be explicitly written.

 \noindent
 In the situations with non-Fredholm operators, the usual solvability
 conditions may not be applicable and the solvability relations
 are, in general, not known. There are some classes of operators
 for which the solvability conditions are obtained. We illustrate
 that with the following example. Consider the equation
 \begin{equation}
 \label{int1}
 Lu \equiv \Delta u + a u = f
 \end{equation}
 in $\mathbb R^d$, where $a$ is a positive constant.
 Clearly, the operator $L$ coincides with its limiting operators. The
 homogeneous problem admits a nonzero bounded solution. Thus the
 Fredholm property is not satisfied. However, since the operator
 has constant coefficients, we can apply the Fourier transform and
 find the solution explicitly. Solvability relations can be
 formulated as follows. If $f \in L^2(\mathbb R^d)$ and 
 $xf \in L^1(\mathbb R^d)$, then there
 exists a unique solution of this problem in $H^2(\mathbb R^d)$ if and
 only if
 $$  
 \Bigg(f(x),\frac{e^{ipx}}{(2\pi)^{\frac{d}{2}}}\Bigg)_{L^2(\mathbb R^d)}=0, \quad 
 p\in S_{\sqrt{a}}^{d} \quad a.e.
 $$
 (see \cite{VV103}). Here $S_{\sqrt{a}}^{d}$ denotes the
 sphere in $\mathbb R^d$ of radius $\sqrt{a}$ centered at the origin.
 Thus, though the operator fails to satisfy the Fredholm property,
 the solvability conditions are formulated similarly. However,
 this similarity is only formal since the range of the operator is
 not closed.

 \noindent
 In the case of the operator with a potential,
 $$
 L u \equiv \Delta u + a(x) u = f ,
 $$
 the Fourier transform is not directly applicable. Nevertheless, the solvability
 conditions in ${\mathbb R}^{3}$ can be obtained by a rather sophisticated 
 application of the theory of self-adjoint
 operators (see \cite{VV08}). As before, the solvability relations are
 formulated 
 in terms of the orthogonality to the solutions of the homogeneous adjoint
 problem. There are several other examples of the linear elliptic
 non-Fredholm operators for which the solvability
 conditions can be obtained (see ~\cite{EV22}, ~\cite{V11},
 ~\cite{VKMP02}, ~\cite{VV08}, ~\cite{VV10}, ~\cite{VV103}).

 \noindent
 Solvability conditions play a crucial role in the analysis of the
 nonlinear elliptic equations. In the case of the operators without Fredholm
 property, in spite of some progress in the understanding of the linear
 problems, there exist only few examples where nonlinear non-Fredholm
 operators are analyzed (see ~\cite{DMV05}, ~\cite{DMV08}, ~\cite{EV20},
 ~\cite{EV22}, ~\cite{VV103}).

 \noindent
 In the present article we study another class of stationary
 nonlinear systems of equations, for which the Fredholm property may or may
 not be satisfied:
\begin{equation}
\label{id1}
\frac{d^{2}u_{k}}{dx^{2}}+ b_{k}\frac{du_{k}}{dx}+a_{k}u_{k} + 
\int_{\Omega}G_{k}(x-y)F_{k}(u_{1}(y), u_{2}(y), ..., u_{N}(y), y)dy = 0,
\quad 1\leq k\leq K,
\end{equation}
\begin{equation}
\label{id2}
\frac{d^{2}u_{k}}{dx^{2}}+ a_{k}u_{k} + 
\int_{\Omega}G_{k}(x-y)F_{k}(u_{1}(y), u_{2}(y), ..., u_{N}(y), y)dy = 0,
\quad K+1\leq k\leq N,
\end{equation}
with the constants
$a_{k}\geq 0, \ b_{k}\in {\mathbb R}, \ b_{k}\neq 0$ and $K\geq 2$.
Here $x\in \Omega\subseteq {\mathbb R}$. In the first case we consider
the situation when $\Omega={\mathbb R}$, such that $N\geq 4$. In the second
part of the article we discuss the case of the finite interval
$\Omega=I:=[0, 2\pi]$ with periodic boundary conditions, so that $N\geq 5$. 
Throughout the work the vector function
\begin{equation}
\label{u}
u:=(u_{1}, u_{2}, ..., u_{N})^{T}\in {\mathbb R}^{N}.
\end{equation}  
For the simplicity of the presentation we restrict ourselves to the one
dimensional
case (the multidimensional case will be considered in our forthcoming paper).
Article ~\cite{EV20} is devoted
to the studies of a single integro-differential equation with a drift term
and the case without a transport term was covered in ~\cite{VV111}.
In the population dynamics the
integro-differential equations describe the models with the intra-specific
competition and the nonlocal consumption of resources (see e.g.
~\cite{ABVV10}, ~\cite{BNPR09}, ~\cite{GVA06}). The studies of
the systems of integro-
differential equations are of interest to us in the context of the complicated
biological systems, where $u_{k}(x,t), \ k=1,...,N$ stand for the cell densities
for various groups of cells in the organism. Let us use the explicit form of
the solvability relations and study the existence of solutions of such
nonlinear systems. We would like to emphasize especially that the solutions of
the integro-differential
equations with the drift terms are relevant to the understanding of the
emergence and propagation of patterns in the theory of speciation
(see ~\cite{VV13}).
The solvability of the linear problems involving the Laplace operator with
the drift term was treated in ~\cite{VV10}, see also ~\cite{BHN05}.
Standing lattice solitons in the discrete NLS equation with saturation
were discussed in ~\cite{AKLP19}. Fredholm structures, topological invariants
and applications were covered in ~\cite{E09}. The work ~\cite{E10} deals with
the finite and infinite dimensional attractors for the evolution equations of
mathematical physics. The large time behavior of solutions of fourth order
parabolic equations and $\varepsilon$-entropy of their attractors were
analyzed in ~\cite{EP07}. The articles ~\cite{GS05} and ~\cite{RS01} are
crucial for the understanding of the Fredholm and properness properties of
the quasilinear elliptic systems of the second order and of the operators
of this kind on ${\mathbb R}^{N}$. The work ~\cite{GS10} is devoted to the
exponential decay and Fredholm properties in the second-order quasilinear
elliptic systems of equations.


\setcounter{equation}{0}

\section{Formulation of the results}

Our technical conditions are analogous to the ones of the ~\cite{EV20},
adapted to the work with vector functions. It is also more difficult to
work in the Sobolev space $H^{2}(\Omega, {\mathbb R}^{N})$, especially when
$\Omega$ is a finite interval with periodic boundary conditions with the
constraints applied.
The nonlinear parts of system (\ref{id1}), (\ref{id2}) will satisfy the
following regularity conditions.

\bigskip

\noindent
{\bf Assumption 2.1.} {\it Let $1\leq k\leq N$. Functions
$F_{k}(u,x): {\mathbb R}^{N}\times \Omega
\to {\mathbb R}$ are satisfying the Caratheodory condition 
(see ~\cite{K64}), such that
\begin{equation}
\label{ub1}
\sqrt{\sum_{k=1}^{N}F_{k}^{2}(u,x)}\leq {\cal K}|u|_{{\mathbb R}^{N}}+h(x) \quad for
\quad u\in {\mathbb R}^{N}, \ x\in \Omega
\end{equation}
with a constant ${\cal K}>0$ and $h(x):\Omega\to {\mathbb R}^{+}, \quad
h(x)\in L^{2}(\Omega)$. Furthermore, they are Lipschitz continuous
functions, such that for any $u^{(1), (2)}\in {\mathbb R}^{N}, \ x\in \Omega:$
\begin{equation}
\label{lk1}
\sqrt{\sum_{k=1}^{N}(F_{k}(u^{(1)},x)-F_{k}(u^{(2)},x))^{2}}\leq L |u^{(1)}-u^{(2)}|_
{{\mathbb R}^{N}},
\end{equation}
with a constant $L>0$.

\noindent
In the case of $\Omega=I$ we assume that $F_{k}(u,0)=F_{k}(u,2\pi)$ for
$u\in {\mathbb R}^{N}$ and all $1\leq k\leq N$.}

\bigskip

\noindent
Here and further down the norm of a vector function given by (\ref{u}) is
$$
|u|_{{\mathbb R}^{N}}:=\sqrt{\sum_{k=1}^{N}u_{k}^{2}}.
$$
Note that the solvability of a local elliptic equation in a bounded domain in
${\mathbb R}^{N}$ was considered in ~\cite{BO86}, where the nonlinear function
was allowed to have a sublinear growth.

\noindent
In order to study the existence of solutions of (\ref{id1}), (\ref{id2}), we
introduce the auxiliary system of equations as
\begin{equation}
\label{ae1}
-\frac{d^{2}u_{k}}{dx^{2}}-b_{k}\frac{du_{k}}{dx}-a_{k}u_{k}=\int_{\Omega}
G_{k}(x-y)F_{k}(v_{1}(y), v_{2}(y), ..., v_{N}(y), y)dy, \quad 1\leq k\leq K,
\end{equation}
\begin{equation}
\label{ae2}
-\frac{d^{2}u_{k}}{dx^{2}}-a_{k}u_{k}=\int_{\Omega}
G_{k}(x-y)F_{k}(v_{1}(y), v_{2}(y), ..., v_{N}(y), y)dy, \quad K+1\leq k\leq N,
\end{equation}
where $a_{k}\geq 0, \ b_{k}\in {\mathbb R}, \ b_{k}\neq 0$ are the constants
and $K\geq 2$.

\noindent
Let us designate
\begin{equation}
\label{ip}  
(f_{1}(x),f_{2}(x))_{L^{2}(\Omega)}:=\int_{\Omega}f_{1}(x)\bar{f_{2}}(x)dx,
\end{equation}
with a slight abuse of notations when these functions are not square
integrable, like for instance those involved in orthogonality relation
(\ref{or1}) below. Indeed, if $f_{1}(x)\in L^{1}(\Omega)$ and $f_{2}(x)$ is
bounded, the integral in the right side of (\ref{ip}) makes sense.

\noindent
Let us first consider the situation of the whole real line, such that
$\Omega={\mathbb R}$.
The appropriate Sobolev space is equipped with the norm
\begin{equation}
\label{h2n}  
\|\phi\|_{H^{2}({\mathbb R})}^{2}:=\|\phi\|_{L^{2}({\mathbb R})}^{2}+
\Bigg\|\frac{d^{2}\phi}{dx^{2}}\Bigg\|_{L^{2}({\mathbb R})}^{2}.
\end{equation}
For a vector function with real valued components given by (\ref{u}), we have
\begin{equation}
\label{h2nv}
\|u\| _{H^{2}({\mathbb R}, {\mathbb R}^{N})}^{2}:=\sum_{k=1}^{N}\|u_{k}\|_{H^{2}({\mathbb R})}^{2}
=\sum_{k=1}^{N}\Bigg\{\|u_{k}\|_{L^{2}({\mathbb R})}^{2}+\Bigg\|\frac{d^{2}u_{k}}{dx^{2}}
\Bigg\|_{L^{2}({\mathbb R})}^{2}\Bigg\}.
\end{equation}
We also use the norm
$$
\|u\| _{L^{2}({\mathbb R}, {\mathbb R}^{N})}^{2}:=\sum_{k=1}^{N}\|u_{k}\|_{L^{2}({\mathbb R})}^{2}.
$$
By means of Assumption 2.1 above, we are not allowed to consider the higher
powers
of the nonlinearities, than the first one, which is restrictive from the point
of view of the applications. But this guarantees that our nonlinear vector
function is a bounded and continuous map from
$L^{2}(\Omega, {\mathbb R}^{N})$ to $L^{2}(\Omega, {\mathbb R}^{N})$.

\noindent
The main issue for the system above is that in the absence
of the drift terms we are dealing with the self-adjoint, non-Fredholm
operators
$$
-\frac{d^{2}}{dx^{2}}-a_{k}: H^{2}({\mathbb R})\to L^{2}({\mathbb R}), \
a_{k}\geq 0,
$$
which is the obstacle to solve our problem. The similar situations but
in linear problems, both self- adjoint and non-self-adjoint
involving the differential operators without the Fredholm property
have been studied extensively in recent years (see ~\cite{EV22},
~\cite{V11}, ~\cite{VKMP02}, ~\cite{VV08}, ~\cite{VV10}, ~\cite{VV103}).
However, the situation differs when
the constants in the drift terms $b_{k}\neq 0$. For $1\leq k\leq K$, the
operators
\begin{equation}
\label{lab}
L_{a, \ b, \ k}:=-\frac{d^{2}}{dx^{2}}-b_{k}\frac{d}{dx}-a_{k}: \quad  H^{2}
({\mathbb R})\to L^{2}({\mathbb R})
\end{equation}
with $a_{k}\geq 0$ and $b_{k}\in {\mathbb R}, \ b_{k}\neq 0$ involved in the
left side of (\ref{ae1}) are non-self-adjoint.

\noindent
By means of the standard
Fourier transform, it can be easily verified that the essential spectra of
the operators $L_{a, \ b, \ k}$ are given by
$$
\lambda_{a, \ b, \ k}(p)=p^{2}-a_{k}-ib_{k}p, \quad p\in {\mathbb R}.
$$
Obviously, when $a_{k}>0$ the operators $L_{a, \ b, \ k}$ are Fredholm, because
their essential spectra stay away from the origin. But when $a_{k}=0$ our
operators
$L_{a, \ b, \ k}$ fail to satisfy the Fredholm property since the origin belongs
to their essential spectra.

\noindent
We manage to establish that under the reasonable
technical assumptions system (\ref{ae1}),  (\ref{ae2}) defines a map
$T_{a, b}: H^{2}({\mathbb R}, {\mathbb R}^{N})\to H^{2}({\mathbb R},
{\mathbb R}^{N})$, which is a strict contraction.

\bigskip

\noindent
{\bf Theorem 2.2.}  {\it Let $\Omega={\mathbb R}, \ N\geq 4, \ K\geq 2, \
1\leq l \leq K-1, \ K+1\leq r\leq N-1$, the integral kernels
$G_{k}(x): {\mathbb R}\to {\mathbb R}, \ G_{k}(x)\in L^{1}({\mathbb R})$
for all $1\leq k\leq N$ and
Assumption 2.1  holds.

\medskip

\noindent
a) Let $a_{k}>0, \ b_{k}\in {\mathbb R}, \ b_{k}\neq 0$ for $1\leq k\leq l$.

\medskip

\noindent
b) Let $a_{k}=0, \ b_{k}\in {\mathbb R}, \ b_{k}\neq 0$ for $l+1\leq k\leq K$,
additionally $xG_{k}(x)\in L^{1}({\mathbb R})$ and orthogonality conditions
(\ref{or1}) hold.

\medskip

\noindent
c) Let $a_{k}>0, \ xG_{k}(x)\in L^{1}({\mathbb R})$ for $K+1\leq k\leq r$ and
orthogonality relations (\ref{or12}) hold.

\medskip

\noindent
d) Let $a_{k}=0, \ x^{2}G_{k}(x)\in L^{1}({\mathbb R})$ for $r+1\leq k\leq N$,
orthogonality conditions (\ref{or13}) hold
and $2\sqrt{\pi}QL<1$ with $Q$ defined in (\ref{Q}) below.
Then the map $v \mapsto T_{a, b}v=u$ on $H^{2}({\mathbb R}, {\mathbb R}^{N})$
defined by problem (\ref{ae1}), (\ref{ae2}) has a unique fixed point
$v^{(a, b)}$, which is
the only solution of the system of equations (\ref{id1}), (\ref{id2}) in
$H^{2}({\mathbb R}, {\mathbb R}^{N})$.

\medskip

\noindent
The fixed point $v^{(a, b)}$ is nontrivial provided that for some $1\leq k\leq N$
the intersection of supports of the Fourier transforms
of functions $\hbox{supp}\widehat{F_{k}(0,x)}\cap \hbox{supp} \widehat{G_{k}}$
is a set of nonzero Lebesgue measure in ${\mathbb R}$.}

\bigskip

\noindent
Let us note that in the case a) of the theorem above, when
$a_{k}>0, \ b_{k}\in {\mathbb R}, \ b_{k}\neq 0$, the orthogonality
relations are not needed.

\noindent
In the second part of the article we study the analogous system on the 
finite interval $I=[0, 2\pi]$ with periodic boundary conditions
with
$a_{k}\geq 0, \ b_{k}\in {\mathbb R}, \ b_{k}\neq 0, \
K\geq 2, \ N\geq 5$, namely for $1\leq k\leq K$
\begin{equation}
\label{id1m}
\frac{d^{2}u_{k}}{dx^{2}}+ b_{k}\frac{du_{k}}{dx}+a_{k}u_{k}+ 
\int_{0}^{2\pi}G_{k}(x-y)F_{k}
(u_{1}(y), u_{2}(y), ...,  u_{N}(y), y)dy = 0,
\end{equation}
and for $K+1\leq k\leq N$
\begin{equation}
\label{id2m}
\frac{d^{2}u_{k}}{dx^{2}}+a_{k}u_{k}+\int_{0}^{2\pi}G_{k}(x-y)F_{k}
(u_{1}(y), u_{2}(y), ...,  u_{N}(y), y)dy = 0.
\end{equation}
Let us use the function space
\begin{equation}
\label{h2i}
H^{2}(I):=\{v(x): I\to {\mathbb R} \ | \ v(x), v''(x)\in L^{2}(I), \
v(0)=v(2\pi), \ v'(0)=v'(2\pi)\}
\end {equation}  
aiming at $u_{k}(x)\in H^{2}(I)$ for $1\leq k\leq l$ and $K+1\leq k\leq r$
with $1\leq l\leq K-1$ and $K+1\leq r\leq q-1$ (see Theorem 2.3 and
Lemma A2 below).

\noindent
For the technical purposes we introduce the auxiliary constrained subspaces
\begin{equation}
\label{h2ik}
H_{k}^{2}(I):= \Big\{v\in H^{2}(I) \ \Big| \ \Big(v(x),
\frac{e^{\pm in_{k}x}}{\sqrt{2\pi}}\Big)_{L^{2}(I)}=0 \Big\}, \quad
n_{k}\in {\mathbb N}, \quad r+1\leq k\leq q.
\end {equation}
Our goal is to have $u_{k}(x)\in H_{k}^{2}(I), \ r+1\leq k\leq q$, where
$r+1\leq q\leq N-1$. Also,
\begin{equation}
\label{h02i}
H_{0}^{2}(I):=\{v\in H^{2}(I) \ | \ (v(x), 1)_{L^{2}(I)}=0\}.
\end{equation}
We plan to have $u_{k}(x)\in H_{0}^{2}(I)$ for $l+1\leq k\leq K$ and
$q+1\leq k\leq N$ (see Theorem 2.3 and Lemma A2). The constrained subspaces
(\ref{h2ik}) and (\ref{h02i}) are Hilbert spaces as well (see e.g. Chapter 2.1
of ~\cite{HS96}).

\noindent
The resulting space used for demonstrating the existence of the solution
$u(x): I \to {\mathbb R}^{N}$ of the system of equations
(\ref{id1m}), (\ref{id2m}) will
be the direct sum of the spaces introduced above, namely
\begin{equation}
\label{hc2i}
H_{c}^{2}(I, {\mathbb R}^{N}):=\oplus_{k=1}^{l}H^{2}(I)\oplus_{k=l+1}^{K}H_{0}^{2}(I)
\oplus_{k=K+1}^{r}H^{2}(I)\oplus_{k=r+1}^{q}H_{k}^{2}(I)\oplus_{k=q+1}^{N}H_{0}^{2}(I).
\end{equation}
The corresponding Sobolev norm will be
\begin{equation}
\label{hc2in}
\|u\|_{H_{c}^{2}(I, {\mathbb R}^{N})}^{2}:=\sum_{k=1}^{N}\{\|u_{k}\|_{L^{2}(I)}^{2}+
\|u_{k}''\|_{L^{2}(I)}^{2}\},
\end{equation}  
where $u(x): I \to {\mathbb R}^{N}$.

\noindent
Let us demonstrate that under the reasonable
technical conditions system (\ref{ae1}),  (\ref{ae2}) with $\Omega=I$
defines a map
$\tau_{a, b}: H_{c}^{2}(I, {\mathbb R}^{N})\to H_{c}^{2}(I, {\mathbb R}^{N})$, which
is a strict contraction.

\bigskip

\noindent
{\bf Theorem 2.3.}  {\it Let $\Omega=I, \ N\geq 5, \ K\geq 2, \
1\leq l\leq K-1, \ K+1\leq r\leq q-1, \ r+1\leq q\leq N-1$, the
integral kernels
$G_{k}(x): I\to {\mathbb R}, \ G_{k}(x)\in C(I), \ G_{k}(0)=G_{k}(2\pi)$ for
all $1\leq k\leq N$ and Assumption 2.1 is valid.

\medskip

\noindent
a) Let $a_{k}>0, \ b_{k}\in {\mathbb R}, \ b_{k}\neq 0$ for $1\leq k\leq l$.

\medskip

\noindent
b) Let $a_{k}=0, \ b_{k}\in {\mathbb R}, \ b_{k}\neq 0$ for $l+1\leq k\leq K$
and orthogonality condition (\ref{or2}) holds.

\medskip

\noindent
c) Let $a_{k}>0, \ a_{k}\neq n^{2}, \ n\in {\mathbb Z}$ for $K+1\leq k\leq r$.

\medskip

\noindent
d) Let $a_{k}=n_{k}^{2}, \ n_{k}\in {\mathbb N}$ for $r+1\leq k\leq q$
and orthogonality relations (\ref{or21}) are valid.

\medskip

\noindent
e) Let $a_{k}=0$ for $q+1\leq k\leq N$, orthogonality condition (\ref{or2})
holds and $2\sqrt{\pi}{\cal Q}L<1$, where ${\cal Q}$ is introduced in
(\ref{Q2}).
Then the map $v \mapsto \tau_{a, b}v=u$ on
$H_{c}^{2}(I, {\mathbb R}^{N})$
defined by the system of equations (\ref{ae1}), (\ref{ae2}) has a unique fixed
point $v^{(a, b)}$, which is
the only solution of system (\ref{id1m}), (\ref{id2m}) in
$H_{c}^{2}(I, {\mathbb R}^{N})$.

\medskip

\noindent
The fixed point $v^{(a, b)}$ does not vanish identically in $I$
provided that for some $1\leq k\leq N$ and a certain $n\in {\mathbb Z}$ the
Fourier coefficients $G_{k, n}F_{k}(0, x)_{n}\neq 0$.}

\bigskip

\noindent
{\bf Remark 2.4.} {\it Note that in the present work we deal with real
valued vector functions by means of the assumptions on
$F_{k}(u,x)$ and $G_{k}(x)$
involved in the nonlocal terms of the systems of
equations considered above.}


\setcounter{equation}{0}

\section{The existence of solutions for the integro-differential systems}

\bigskip

\noindent
{\it Proof of Theorem 2.2.} First we suppose that 
for a certain $v\in H^{2}({\mathbb R}, {\mathbb R}^{N})$
there exist two solutions $u^{(1),(2)}\in H^{2}({\mathbb R}, {\mathbb R}^{N})$ of
problem (\ref{ae1}), (\ref{ae2}). Then their difference
$w(x):=u^{(1)}(x)-u^{(2)}(x)\in  H^{2}({\mathbb R}, {\mathbb R}^{N})$ will be
a solution of the homogeneous system of equations
$$
-\frac{d^{2}w_{k}}{dx^{2}}-b_{k}\frac{dw_{k}}{dx}-a_{k}w_{k}=0, \quad 1\leq k\leq K,
$$
$$
-\frac{d^{2}w_{k}}{dx^{2}}-a_{k}w_{k}=0, \quad K+1\leq k\leq N.
$$
But the operators
$\displaystyle{-\frac{d^{2}}{dx^{2}}-a_{k}, \ L_{a, \ b, \ k}: H^{2}({\mathbb R})\to
L^{2}({\mathbb R})}$, where  $L_{a, \ b, \ k}$ is defined in (\ref{lab}) do not
have any nontrivial zero modes. Therefore,
$w(x)$ vanishes on  ${\mathbb R}$.

\noindent
Let us choose an arbitrary $v(x)\in H^{2}({\mathbb R}, {\mathbb R}^{N})$. We
apply the standard Fourier transform (\ref{ft}) to both sides of system
(\ref{ae1}), (\ref{ae2}). This gives us
\begin{equation}
\label{f1}
\widehat{u_{k}}(p)=\sqrt{2\pi}{\widehat{G_{k}}(p)\widehat{f_{k}}(p)\over
p^{2}-a_{k}-ib_{k}p}, \quad
p^{2}\widehat{u_{k}}(p)=\sqrt{2\pi}{p^{2}\widehat{G_{k}}(p)\widehat{f_{k}}(p)\over
p^{2}-a_{k}-ib_{k}p}, \quad 1\leq k\leq K,
\end{equation}
\begin{equation}
\label{f11}
\widehat{u_{k}}(p)=\sqrt{2\pi}{\widehat{G_{k}}(p)\widehat{f_{k}}(p)\over
p^{2}-a_{k}}, \quad
p^{2}\widehat{u_{k}}(p)=\sqrt{2\pi}{p^{2}\widehat{G_{k}}(p)\widehat{f_{k}}(p)\over
p^{2}-a_{k}}, \quad K+1\leq k\leq N.
\end{equation}
Here $\widehat{f_{k}}(p)$ stands for the Fourier image of $F_{k}(v(x),x)$.

\noindent
Clearly, for $1\leq k\leq K$, we have the upper bounds with
$N_{a, \ b, \ k}$ defined in (\ref{Na}), namely
$$
|\widehat{u_{k}}(p)|\leq \sqrt{2\pi}N_{a, \ b, \ k}|\widehat{f_{k}}(p)|  \quad and
\quad |p^{2}\widehat{u_{k}}(p)|\leq \sqrt{2\pi}N_{a, \ b, \ k}|\widehat{f_{k}}(p)|.
$$
For $K+1\leq k\leq N$, we obtain 
$$
|\widehat{u_{k}}(p)|\leq \sqrt{2\pi}M_{a, \ k}|\widehat{f_{k}}(p)|  \quad and
\quad |p^{2}\widehat{u_{k}}(p)|\leq \sqrt{2\pi}M_{a, \ k}|\widehat{f_{k}}(p)|,
$$
with $M_{a, \ k}$ introduced in (\ref{Ma}).
Note that $N_{a, \ b, \ k}<\infty$ by virtue of Lemma A1 of the Appendix without
any orthogonality conditions for $a_{k}>0, \ 1\leq k\leq l$ and under
orthogonality relation
(\ref{or1}) when $a_{k}=0, \ l+1\leq k\leq K$. Also, $M_{a, \ k}<\infty$ under
orthogonality
conditions (\ref{or12}) when $a_{k}>0, \ K+1\leq k\leq r$ and under
orthogonality relations
(\ref{or13}) for $a_{k}=0, \ r+1\leq k\leq N$.

\noindent
This allows us to obtain the upper bound on the norm as
$$
\|u\|_{H^{2}({\mathbb R}, {\mathbb R}^{N})}^{2}=\sum_{k=1}^{N}\{
\|\widehat{u_{k}}(p)\|_{L^{2}({\mathbb R})}^{2}+
\|p^{2}\widehat{u_{k}}(p)\|_{L^{2}({\mathbb R})}^{2}\}\leq
$$
\begin{equation}
\label{nmabk}  
\leq 4\pi\sum_{k=1}^{K}
N_{a, \ b, \ k}^{2}\|F_{k}(v(x),x)\|_{L^{2}({\mathbb R})}^{2}+
4\pi\sum_{k=K+1}^{N}
M_{a, \ k}^{2}\|F_{k}(v(x),x)\|_{L^{2}({\mathbb R})}^{2}.
\end{equation}
The right side of (\ref{nmabk}) is finite via (\ref{ub1}) of Assumption 2.1
since $|v(x)|_{{\mathbb R}^{N}}\in L^{2}({\mathbb R})$. Thus, for an arbitrary
$v(x)\in H^{2}({\mathbb R}, {\mathbb R}^{N})$ there exists a unique solution
$u(x)\in H^{2}({\mathbb R}, {\mathbb R}^{N})$ of problem (\ref{ae1}),
(\ref{ae2}).
Its Fourier image is given by (\ref{f1}), (\ref{f11}). Therefore, the map
$T_{a, b}:  H^{2}({\mathbb R}, {\mathbb R}^{N})\to  H^{2}({\mathbb R},
{\mathbb R}^{N})$ is well defined.

\noindent
This enables us to choose arbitrarily
$v^{(1), (2)}(x)\in H^{2}({\mathbb R}, {\mathbb R}^{N})$,
such that under the given conditions their images
$u^{(1), (2)}:=T_{a, b}v^{(1),(2)}\in H^{2}({\mathbb R}, {\mathbb R}^{N})$.
By means of (\ref{ae1}), (\ref{ae2}) along with (\ref{ft}),
$$
\widehat{u_{k}^{(1)}}(p)=\sqrt{2\pi}{\widehat{G_{k}}(p)\widehat{f_{k}^{(1)}}(p)\over
p^{2}-a_{k}-ib_{k}p}, \quad
\widehat{u_{k}^{(2)}}(p)=\sqrt{2\pi}{\widehat{G_{k}}(p)\widehat{f_{k}^{(2)}}(p)\over
p^{2}-a_{k}-ib_{k}p}, \quad
1\leq k\leq K,
$$
$$
\widehat{u_{k}^{(1)}}(p)=\sqrt{2\pi}{\widehat{G_{k}}(p)\widehat{f_{k}^{(1)}}(p)\over
p^{2}-a_{k}}, \quad
\widehat{u_{k}^{(2)}}(p)=\sqrt{2\pi}{\widehat{G_{k}}(p)\widehat{f_{k}^{(2)}}(p)\over
p^{2}-a_{k}}, \quad K+1\leq k\leq N.
$$
Here $\widehat{f_{k}^{(1)}}(p)$ and $\widehat{f_{k}^{(2)}}(p)$  denote the
Fourier transforms of
$F_{k}(v^{(1)}(x),x)$  and $F_{k}(v^{(2)}(x),x)$ respectively.

\noindent
Hence, for
$1\leq k\leq K$, we easily derive
$$
\Bigg|\widehat{u_{k}^{(1)}}(p)-\widehat{u_{k}^{(2)}}(p)\Bigg|\leq
\sqrt{2\pi}N_{a, \ b, \ k}
\Bigg|\widehat{f_{k}^{(1)}}(p)-\widehat{f_{k}^{(2)}}(p)\Bigg|,
$$
$$
\Bigg|p^{2}\widehat{u_{k}^{(1)}}(p)-p^{2}\widehat{u_{k}^{(2)}}(p)\Bigg|\leq
\sqrt{2\pi}N_{a, \ b, \ k}
\Bigg|\widehat{f_{k}^{(1)}}(p)-\widehat{f_{k}^{(2)}}(p)\Bigg|
$$
and for $K+1\leq k\leq N$
$$
\Bigg|\widehat{u_{k}^{(1)}}(p)-\widehat{u_{k}^{(2)}}(p)\Bigg|\leq
\sqrt{2\pi}M_{a, \ k}
\Bigg|\widehat{f_{k}^{(1)}}(p)-\widehat{f_{k}^{(2)}}(p)\Bigg|,
$$
$$
\Bigg|p^{2}\widehat{u_{k}^{(1)}}(p)-p^{2}\widehat{u_{k}^{(2)}}(p)\Bigg|\leq
\sqrt{2\pi}M_{a, \ k}
\Bigg|\widehat{f_{k}^{(1)}}(p)-\widehat{f_{k}^{(2)}}(p)\Bigg|.
$$
Then for the appropriate norm of the difference of vector functions we obtain
$$
\|u^{(1)}-u^{(2)}\|_{H^{2}({\mathbb R}, {\mathbb R}^{N})}^{2}=\sum_{k=1}^{N}\Big\{
\Big\|\widehat{u_{k}^{(1)}}(p)-\widehat{u_{k}^{(2)}}(p)\Big\|_{L^{2}({\mathbb R})}^{2}+
\Big\|p^{2}\widehat{u_{k}^{(1)}}(p)-p^{2}\widehat{u_{k}^{(2)}}(p)\Big\|
_{L^{2}({\mathbb R})}^{2}\Big\}\leq
$$
$$
\leq 4\pi Q^{2}
\sum_{k=1}^{N}\|F_{k}(v^{(1)}(x),x)-F_{k}(v^{(2)}(x),x)\|_{L^{2}({\mathbb R})}^{2},
$$
where $Q$ is defined in (\ref{Q}).
Clearly, all $v_{k}^{(1), (2)}(x)\in H^{2}({\mathbb R})\subset
L^{\infty}({\mathbb R})$ due to the Sobolev embedding.

\noindent
Condition (\ref{lk1}) gives us
$$
\sum_{k=1}^{N}\|F_{k}(v^{(1)}(x),x)-F_{k}(v^{(2)}(x),x)\|_{L^{2}({\mathbb R})}^{2}\leq
L^{2}\|v^{(1)}-v^{(2)}\|_{L^{2}({\mathbb R}, {\mathbb R}^{N})}^{2}.
$$
Thus,
\begin{equation}
\label{QL}  
\|T_{a, b}v^{(1)}-T_{a, b}v^{(2)}\|_{H^{2}({\mathbb R}, {\mathbb R}^{N})}\leq 2
\sqrt{\pi}QL\|v^{(1)}-v^{(2)}\|_{H^{2}({\mathbb R}, {\mathbb R}^{N})}.
\end{equation}
The constant in the right side of (\ref{QL}) is less than one via the one of 
our assumptions. Therefore, by virtue of the Fixed Point Theorem,
there exists a unique vector function
$v^{(a, b)}\in H^{2}({\mathbb R}, {\mathbb R}^{N})$, such that
$T_{a, b}v^{(a, b)}=v^{(a, b)}$,
which is the only solution of problem (\ref{id1}), (\ref{id2}) in
$H^{2}({\mathbb R}, {\mathbb R}^{N})$. Suppose $v^{(a, b)}(x)$ vanishes
identically on the real line.
This will contradict to our assumption that for a certain $1\leq k\leq N$, the
Fourier transforms of $G_{k}(x)$ and $F_{k}(0,x)$ are nontrivial on a set
of nonzero Lebesgue measure in ${\mathbb R}$.
\hfill\lanbox

\bigskip

\noindent
Then we turn our attention to establishing the existence of
solutions for our system of integro-differential equations on the finite
interval with periodic boundary conditions.

\bigskip

\noindent
{\it Proof of Theorem 2.3.} Let us first suppose that for some
$v\in H_{c}^{2}(I, {\mathbb R}^{N})$ there exist two solutions
$u^{(1), (2)}\in H_{c}^{2}(I, {\mathbb R}^{N})$ of system (\ref{ae1}), (\ref{ae2})
with $\Omega=I$. Then the difference
${\tilde w}(x):=u^{(1)}(x)-u^{(2)}(x)\in H_{c}^{2}(I, {\mathbb R}^{N})$ will
satisfy the homogeneous system of equations
\begin{equation}
\label{hsi1}  
-\frac{d^{2}{\tilde w}_{k}}{dx^{2}}-b_{k}\frac{d{\tilde w}_{k}}{dx}-
a_{k}{\tilde w}_{k}=0, \quad 1\leq k\leq K,
\end{equation}
\begin{equation}
\label{hsi2}  
-\frac{d^{2}{\tilde w}_{k}}{dx^{2}}-a_{k}{\tilde w}_{k}=0, \quad K+1\leq k\leq N.
\end{equation}
Evidently, each operator contained in the left side of system (\ref{hsi1})
\begin{equation}
\label{clabk}
{\cal L}_{a, \ b, \ k}:=-\frac{d^{2}}{dx^{2}}-b_{k}\frac{d}{dx}-a_{k}: \quad
H^{2}(I)\to L^{2}(I),      
\end{equation}
where $1\leq k\leq l, \ a_{k}>0, \ b_{k}\in {\mathbb R}, \ b_{k}\neq 0$ is
Fredholm, non-self-adjoint.

\noindent
Its set of eigenvalues is
\begin{equation}
\label{ev}
\lambda_{a, \ b, \ k}(n)=n^{2}-a_{k}-ib_{k}n, \quad n\in {\mathbb Z}.
\end{equation}
The corresponding eigenfunctions are the standard Fourier harmonics
\begin{equation}
\label{fg}  
\frac{e^{inx}}{\sqrt{2\pi}}, \quad n\in {\mathbb Z}.
\end{equation}
When $l+1\leq k\leq K$, we have
$a_{k}=0, \ b_{k}\in {\mathbb R}, \ b_{k}\neq 0$.
Let us deal with the operators
\begin{equation}
\label{cl0bk}
{\cal L}_{0, \ b, \ k}:=-\frac{d^{2}}{dx^{2}}-b_{k}\frac{d}{dx}: \quad
H_{0}^{2}(I)\to L^{2}(I)     
\end{equation}
involved in system (\ref{hsi1}).

\noindent
Obviously, each operator (\ref{cl0bk})
has the eigenvalues given by formula (\ref{ev}) with
$a_{k}=0, \ b_{k}\in {\mathbb R}, \ b_{k}\neq 0, \ n\in {\mathbb Z}, \ n\neq 0$.

\noindent
The corresponding eigenfunctions in this case are (\ref{fg}) with
$n\in {\mathbb Z}, \ n\neq 0$.

\noindent
Clearly, every operator contained in the left side of the system of equations
(\ref{hsi2})
\begin{equation}
\label{cla0k}
{\cal L}_{a, \ 0, \ k}:=-\frac{d^{2}}{dx^{2}}-a_{k}: \quad H^{2}(I)\to L^{2}(I)      
\end{equation}
with $K+1\leq k\leq r, \ a_{k}>0, \ a_{k}\neq n^{2}, \ n\in {\mathbb Z}$ is
Fredholm, self-adjoint.

\noindent
Its set of eigenvalues is
\begin{equation}
\label{ev0}
\lambda_{a, \ 0, \ k}(n)=n^{2}-a_{k}, \quad n\in {\mathbb Z}.
\end{equation}
The corresponding eigenfunctions are given by (\ref{fg}).

\noindent
For $r+1\leq k\leq q$, we have $a_{k}=n_{k}^{2}, \ n_{k}\in {\mathbb N}$.
Let us consider the operators
\begin{equation}
\label{clnk0k}
{\cal L}_{n_{k}^{2}, \ 0, \ k}:=-\frac{d^{2}}{dx^{2}}-n_{k}^{2}: \quad
H_{k}^{2}(I)\to L^{2}(I)      
\end{equation}
involved in the left side of system (\ref{hsi2}).

\noindent
The eigenvalues of each operator (\ref{clnk0k}) are written in (\ref{ev0})
with $a_{k}=n_{k}^{2}, \ n\in {\mathbb Z}, \ n\neq \pm n_{k}$.

\noindent
The corresponding eigenfunctions are given by formula (\ref{fg}) with
$n\in {\mathbb Z}, \ n\neq \pm n_{k}$.

\noindent
When $q+1\leq k\leq N$, all the constants $a_{k}$ are trivial. Then the operator
contained in the left side of the system of equations (\ref{hsi2}) in this
situation is
\begin{equation}
\label{cl00k}
{\cal L}_{0, \ 0, \ k}:=-\frac{d^{2}}{dx^{2}}: \quad H_{0}^{2}(I)\to L^{2}(I).      
\end{equation}
Its eigenvalues are
\begin{equation}
\label{l00k}
\lambda_{0, \ 0, \ k}(n)=n^{2}, \quad n\in {\mathbb Z}, \quad n\neq 0.  
\end{equation}
We have the corresponding eigenfunctions written in (\ref{fg}) with
$n\in {\mathbb Z}, \ n\neq 0$.

\noindent
Note that all the operators mentioned above, which are involved in the left
side of the homogeneous system (\ref{hsi1}), (\ref{hsi2}) have the trivial
kernels. Thus, the vector function ${\tilde w}(x)$ vanishes identically on
the interval $I$.

\noindent
Let us choose arbitrarily $v(x)\in H_{c}^{2}(I, {\mathbb R}^{N})$ and apply
the Fourier transform (\ref{fti}) to both sides of the system of equations 
(\ref{ae1}), (\ref{ae2}) with $\Omega=I$. This yields
\begin{equation}
\label{ukn1k}
u_{k, n}=\sqrt{2\pi}\frac{G_{k, n}f_{k, n}}{n^{2}-a_{k}-ib_{k}n}, \quad
n^{2}u_{k, n}=\sqrt{2\pi}\frac{n^{2}G_{k, n}f_{k, n}}{n^{2}-a_{k}-ib_{k}n}, \quad
1\leq k\leq K, \quad n\in {\mathbb Z},
\end{equation}
\begin{equation}
\label{uknk1N}
u_{k, n}=\sqrt{2\pi}\frac{G_{k, n}f_{k, n}}{n^{2}-a_{k}}, \quad
n^{2}u_{k, n}=\sqrt{2\pi}\frac{n^{2}G_{k, n}f_{k, n}}{n^{2}-a_{k}}, \quad
K+1\leq k\leq N, \quad n\in {\mathbb Z},
\end{equation}
where $f_{k, n}:=F_{k}(v(x), x)_{n}$.

\noindent
From (\ref{ukn1k}) and (\ref{uknk1N}) we easily obtain that
\begin{equation}
\label{ukn1kub}  
|u_{k, n}|\leq \sqrt{2\pi}{\cal N}_{a, \ b, \ k}|f_{k, n}|, \quad
|n^{2}u_{k, n}|\leq \sqrt{2\pi}{\cal N}_{a, \ b, \ k}|f_{k, n}|, \quad 1\leq k\leq K,
\quad n\in {\mathbb Z},
\end{equation}
\begin{equation}
\label{uknk1kub} 
|u_{k, n}|\leq \sqrt{2\pi}{\cal M}_{a, \ k}|f_{k, n}|, \quad
|n^{2}u_{k, n}|\leq \sqrt{2\pi}{\cal M}_{a,  \ k}|f_{k, n}|, \quad K+1\leq k\leq N,
\quad n\in {\mathbb Z}.
\end{equation}
In the estimates above we have ${\cal N}_{a, \ b, \ k}$ introduced in
(\ref{Nac}) and ${\cal M}_{a,  \ k}$ was defined in (\ref{Mac}). Clearly,
${\cal N}_{a, \ b, \ k}<\infty$ by means of Lemma A2 of the Appendix without
any orthogonality relations for $a_{k}>0, \ 1\leq k\leq l$ and under
orthogonality condition (\ref{or2}) when $a_{k}=0, \ l+1\leq k\leq K$.
Similarly, ${\cal M}_{a,  \ k}<\infty$ for 
$a_{k}>0, \ a_{k}\neq n^{2}, \ n\in {\mathbb Z}, \ K+1\leq k\leq r$,
under orthogonality relations (\ref{or21}) when 
$a_{k}=n_{k}^{2}, \ n_{k}\in {\mathbb N}, \ r+1\leq k\leq q$ and 
under orthogonality condition (\ref{or2})
for $a_{k}=0, \ q+1\leq k\leq N$ by virtue of Lemma A2.

\noindent
By means of (\ref{ukn1kub}) and (\ref{uknk1kub}), we derive
$$
\|u\|_{H_{c}^{2}(I, {\mathbb R}^{N})}^{2}=\sum_{k=1}^{K}\Big[\sum_{n=-\infty}^{\infty}
|u_{k, n}|^{2}+\sum_{n=-\infty}^{\infty}|n^{2}u_{k, n}|^{2}\Big]+
\sum_{k=K+1}^{N}\Big[\sum_{n=-\infty}^{\infty}
|u_{k, n}|^{2}+\sum_{n=-\infty}^{\infty}|n^{2}u_{k, n}|^{2}\Big]\leq
$$
\begin{equation}
\label{nabkmak}
\leq 4\pi \sum_{k=1}^{K}{\cal N}_{a, \ b, \ k}^{2}\|F_{k}(v(x), x)\|_{L^{2}(I)}^{2}+
4\pi \sum_{k=K+1}^{N}{\cal M}_{a, \ k}^{2}\|F_{k}(v(x), x)\|_{L^{2}(I)}^{2}.
\end{equation}
Let us recall inequality (\ref{ub1}) of Assumption 2.1. We have here
$|v(x)|_{{\mathbb R}^{N}}\in L^{2}(I)$, such that all
$F_{k}(v(x), x)\in L^{2}(I)$ and the right side of (\ref{nabkmak}) is finite.
Hence, for any $v(x)\in H_{c}^{2}(I, {\mathbb R}^{N})$ there exists a unique
$u(x)\in H_{c}^{2}(I, {\mathbb R}^{N})$, which solves system
(\ref{ae1}), (\ref{ae2}) with $\Omega=I$. Its Fourier transform is given
by formulas (\ref{ukn1k}) and (\ref{uknk1N}). Thus, the map
$\tau_{a, b}: H_{c}^{2}(I, {\mathbb R}^{N})\to H_{c}^{2}(I, {\mathbb R}^{N})$
is well defined.

\noindent
Let us choose arbitrarily $v^{(1), (2)}(x)\in H_{c}^{2}(I, {\mathbb R}^{N})$.
Under the stated assumptions, their images under the map discussed above
are $u^{(1), (2)}:=\tau_{a, b}v^{(1), (2)}\in H_{c}^{2}(I, {\mathbb R}^{N})$.
By virtue of (\ref{ae1}) and (\ref{ae2}) with $\Omega=I$ along with
(\ref{fti}), we arrive at
$$
u_{k, n}^{(1)}=\sqrt{2\pi}\frac{G_{k, n}f_{k, n}^{(1)}}{n^{2}-a_{k}-ib_{k}n}, \quad
u_{k, n}^{(2)}=\sqrt{2\pi}\frac{G_{k, n}f_{k, n}^{(2)}}{n^{2}-a_{k}-ib_{k}n}, \quad
1\leq k\leq K, \quad n\in {\mathbb Z},
$$
$$
u_{k, n}^{(1)}=\sqrt{2\pi}\frac{G_{k, n}f_{k, n}^{(1)}}{n^{2}-a_{k}}, \quad
u_{k, n}^{(2)}=\sqrt{2\pi}\frac{G_{k, n}f_{k, n}^{(2)}}{n^{2}-a_{k}}, \quad
K+1\leq k\leq N, \quad n\in {\mathbb Z}.
$$
Here $f_{k, n}^{(1)}$ and $f_{k, n}^{(2)}$ stand for the Fourier images of
$F_{k}(v^{(1)}(x), x)$ and $F_{k}(v^{(2)}(x), x)$ respectively under transform
(\ref{fti}).

\noindent
Thus, for $1\leq k\leq K, \ n\in {\mathbb Z}$, we have
$$
|u_{k, n}^{(1)}-u_{k, n}^{(2)}|\leq \sqrt{2\pi}{\cal N}_{a, \ b, \ k}
|f_{k, n}^{(1)}-f_{k, n}^{(2)}|, \quad
|n^{2}u_{k, n}^{(1)}-n^{2}u_{k, n}^{(2)}|\leq \sqrt{2\pi}{\cal N}_{a, \ b, \ k}
|f_{k, n}^{(1)}-f_{k, n}^{(2)}|.
$$
Similarly, for $K+1\leq k\leq N, \ n\in {\mathbb Z}$, we obtain
$$
|u_{k, n}^{(1)}-u_{k, n}^{(2)}|\leq \sqrt{2\pi}{\cal M}_{a, \ k}
|f_{k, n}^{(1)}-f_{k, n}^{(2)}|, \quad
|n^{2}u_{k, n}^{(1)}-n^{2}u_{k, n}^{(2)}|\leq \sqrt{2\pi}{\cal M}_{a, \ k}
|f_{k, n}^{(1)}-f_{k, n}^{(2)}|.
$$
Let us estimate the appropriate norm of the difference of the vector functions
as
$$
\|u^{(1)}-u^{(2)}\|_{H_{c}^{2}(I, {\mathbb R}^{N})}^{2}=\sum_{k=1}^{K}\Big[\sum_{n=-\infty}^
{\infty}|u_{k, n}^{(1)}-u_{k, n}^{(2)}|^{2}+
\sum_{n=-\infty}^{\infty}|n^{2}u_{k, n}^{(1)}-n^{2}u_{k, n}^{(2)}|^{2}\Big]+
$$
$$
+\sum_{k=K+1}^{N}\Big[\sum_{n=-\infty}^
{\infty}|u_{k, n}^{(1)}-u_{k, n}^{(2)}|^{2}+
\sum_{n=-\infty}^{\infty}|n^{2}u_{k, n}^{(1)}-n^{2}u_{k, n}^{(2)}|^{2}\Big]\leq
$$
$$
\leq 4\pi {\cal Q}^{2}\sum_{k=1}^{N}\|F_{k}(v^{(1)}(x), x)-F_{k}(v^{(2)}(x), x)\|_
{L^{2}(I)}^{2},
$$
with ${\cal Q}$ introduced in (\ref{Q2}). Evidently, all
$v_{k}^{(1), (2)}(x)\in H^{2}(I)\subset L^{\infty}(I)$ via the Sobolev embedding.

\noindent
We recall condition (\ref{lk1}) of Assumption 2.1, such that
$$
\sum_{k=1}^{N}\|F_{k}(v^{(1)}(x), x)-F_{k}(v^{(2)}(x), x)\|_{L^{2}(I)}^{2}\leq L^{2}
\|v^{(1)}-v^{(2)}\|_{H_{c}^{2}(I, {\mathbb R}^{N})}^{2}.
$$
Hence,
\begin{equation}
\label{tauabhc2n}  
\|\tau_{a, b}v^{(1)}-\tau_{a, b}v^{(2)}\|_{H_{c}^{2}(I, {\mathbb R}^{N})}\leq
2\sqrt{\pi}{\cal Q}L\|v^{(1)}-v^{(2)}\|_{H_{c}^{2}(I, {\mathbb R}^{N})}.
\end{equation}
The constant in the right side of bound (\ref{tauabhc2n}) is less than one as
assumed. By means of the Fixed Point Theorem, there exists a unique vector
function $v^{(a, b)}\in H_{c}^{2}(I, {\mathbb R}^{N})$, so that
$\tau_{a, b}v^{(a, b)}=v^{(a, b)}$. This is the only solution of the system of
equations (\ref{id1m}), (\ref{id2m}) in $H_{c}^{2}(I, {\mathbb R}^{N})$.
Let us suppose that $v^{(a, b)}(x)$ vanishes identically in I. This will
contradict to the given condition that for a certain $1\leq k\leq N$ and
some $n\in {\mathbb Z}$ the Fourier coefficients
$G_{k, n}F_{k}(0, x)_{n}\neq 0$. \hfill\lanbox


\setcounter{equation}{0}

\section{Appendix}

\bigskip

Let $G_{k}(x)$ be a function,
$G_{k}(x): {\mathbb R}\to {\mathbb R}$,
for which we designate its standard Fourier transform using the hat symbol as
\begin{equation}
\label{ft}  
\widehat{G_{k}}(p):={1\over \sqrt{2\pi}}\int_{-\infty}^{\infty}G_{k}(x)
e^{-ipx}dx, \quad p\in {\mathbb R}.
\end{equation}
Clearly
\begin{equation}
\label{inf1}
\|\widehat{G_{k}}(p)\|_{L^{\infty}({\mathbb R})}\leq {1\over \sqrt{2\pi}}
\|G_{k}\|_{L^{1}({\mathbb R})}
\end{equation}
and
$\displaystyle{G_{k}(x)={1\over \sqrt{2\pi}}\int_{-\infty}^{\infty}
\widehat{G_{k}}(q)e^{iqx}dq, \ x\in {\mathbb R}.}$
For the technical purposes we introduce the auxiliary quantities
\begin{equation}
\label{Na}
N_{a, \ b, \ k}:=\hbox{max}\Big\{
\Big\|{\widehat{G_{k}}(p)\over p^{2}-a_{k}-ib_{k}p}\Big\|_{L^{\infty}({\mathbb R})},
\quad
\Big\|{p^{2}\widehat{G_{k}}(p)\over p^{2}-a_{k}-ib_{k}p}\Big\|_{L^{\infty}({\mathbb R})}
\Big\}, \quad 1\leq k\leq K,
\end{equation}
\begin{equation}
\label{Ma}
M_{a, \ k}:=\hbox{max}\Big\{
\Big\|{\widehat{G_{k}}(p)\over p^{2}-a_{k}}\Big\|_{L^{\infty}({\mathbb R})},
\quad
\Big\|{p^{2}\widehat{G_{k}}(p)\over p^{2}-a_{k}}\Big\|_{L^{\infty}({\mathbb R})}
\Big\}, \quad K+1\leq k\leq N,
\end{equation}
where $a_{k}\geq 0, \ b_{k}\in {\mathbb R}, \ b_{k}\neq 0$ are the constants
and $K\geq 2, \ N \geq 4$.
Let $N_{0, \ b, \ k}$ denote (\ref{Na}) when $a_{k}$ vanishes and $M_{0, \ k}$
stands for (\ref{Ma}) when $a_{k}=0$.
Under the assumptions of Lemma A1 below, quantities (\ref{Na}) and (\ref{Ma})
will be finite. This will allow us to define
\begin{equation}
\label{NaN}
N_{a, \ b}:=\hbox{max}_{1\leq k\leq K}N_{a, \ b, \ k}<\infty,  
\end{equation}
\begin{equation}
\label{MaM}
M_{a}:=\hbox{max}_{K+1\leq k\leq N}M_{a, \ k}<\infty 
\end{equation}
and
\begin{equation}
\label{Q}
Q:=\hbox{max}\{N_{a, \ b}, \ M_{a}\}.
\end{equation}  
The auxiliary lemma below is the adaptation of the ones established in
~\cite{EV20} and ~\cite{VV111} for the studies of the single
integro-differential equation with and without a drift. Let us present it for
the convenience of the readers.
  
\bigskip

\noindent
{\bf Lemma A1.} {\it Let $N\geq 4, \ K\geq 2, \ 1\leq l\leq K-1, \ K+1\leq r\leq
N-1$, the integral kernels       
$G_{k}(x): {\mathbb R}\to {\mathbb R}, \ G_{k}(x)\in L^{1}({\mathbb R})$ for all
$1\leq k\leq N$.

\medskip
  
\noindent  
a) Let $a_{k}>0, \ b_{k}\in {\mathbb R}, \ b_{k}\neq 0$ for $1\leq k\leq l$.
Then $N_{a, \ b, \ k}<\infty$.

\medskip

\noindent
b) Let $a_{k}=0, \ b_{k}\in {\mathbb R}, \ b_{k}\neq 0$ for $l+1\leq k\leq K$
and additionally
$xG_{k}(x)\in L^{1}({\mathbb R})$. Then $N_{0, \ b, \ k}<\infty$ if and only if
\begin{equation}
\label{or1}
(G_{k}(x),1)_{L^{2}({\mathbb R})}=0 
\end{equation}
holds.

\medskip

\noindent
c) Let $a_{k}>0$ and $xG_{k}(x)\in L^{1}({\mathbb R})$ for $K+1\leq k\leq r$.
Then $M_{a, \ k}<\infty$ if and only if
\begin{equation}
\label{or12}
\Bigg(G_{k}(x), \frac{e^{\pm i\sqrt{a_{k}}x}}{\sqrt{2\pi}}\Bigg)_{L^{2}({\mathbb R})}=0  
\end{equation}
is valid.

\medskip

\noindent
d) Let $a_{k}=0$ and $x^{2}G_{k}(x)\in L^{1}({\mathbb R})$ for $r+1\leq k\leq N$.
Then $M_{0, \ k}<\infty$ if and only if
\begin{equation}
\label{or13}
(G_{k}(x), 1)_{L^{2}({\mathbb R})}=0 \quad  and \quad (G_{k}(x), x)_{L^{2}({\mathbb R})}=0
\end{equation}
holds.
}

\medskip

\noindent
{\it Proof.}  Note that in both cases a) and b) of our
lemma the boundedness of
$\displaystyle{\frac{\widehat{G_{k}}(p)} {p^{2}-a_{k}-ib_{k}p}}$ implies the
boundedness of
$\displaystyle{\frac{p^{2}\widehat{G_{k}}(p)}{p^{2}-a_{k}-ib_{k}p}}$. Indeed, we
can express $\displaystyle{\frac{p^{2}\widehat{G_{k}}(p)}{p^{2}-a_{k}-ib_{k}p}}$ as
the following sum
\begin{equation}
\label{sG}
\widehat{G_{k}}(p)+a_{k}\frac{\widehat{G_{k}}(p)}{p^{2}-a_{k}-ib_{k}p}+ib_{k}
\frac{p\widehat{G_{k}}(p)}{p^{2}-a_{k}-ib_{k}p}.
\end{equation}
Obviously, the first term in (\ref{sG}) is bounded by virtue of (\ref{inf1})
because $G_{k}(x)\in L^{1}({\mathbb R})$ due to the one of our assumptions. The
third term in (\ref{sG}) can be trivially estimated from above in the absolute
value by means of (\ref{inf1}) as
$$
\frac{|b_{k}||p||\widehat{G_{k}}(p)|}{\sqrt{(p^{2}-a_{k})^{2}+b_{k}^{2}p^{2}}}\leq
\frac{1}{\sqrt{2\pi}}\|G_{k}(x)\|_{L^{1}({\mathbb R})}<\infty.
$$
Therefore,
$\displaystyle{\frac{\widehat{G_{k}}(p)} {p^{2}-a_{k}-ib_{k}p}\in L^{\infty}
({\mathbb R})}$
yields $\displaystyle{\frac{p^{2}\widehat{G_{k}}(p)} {p^{2}-a_{k}-ib_{k}p}
\in L^{\infty}({\mathbb R})}$.

\noindent
To establish the result of the part a) of our
lemma, we need to estimate
\begin{equation}
\label{Gabp}
\frac{|\widehat{G_{k}}(p)|}{\sqrt{(p^{2}-a_{k})^{2}+b_{k}^{2}p^{2}}}. 
\end{equation}
Evidently, the numerator of (\ref{Gabp}) can be easily estimated from above by
means of (\ref{inf1}) and the denominator in (\ref{Gabp}) can be trivially
bounded below by a finite, positive constant, such that 
$$
\Bigg|\frac{\widehat{G_{k}}(p)} {p^{2}-a_{k}-ib_{k}p}\Bigg|\leq C\|G_{k}(x)\|_
{L^{1}({\mathbb R})}<\infty.
$$
Here and further down $C$ will stand for a finite, positive constant. This
implies that under our assumptions, when $a_{k}>0$ we have $N_{a, \ b, \ k}<\infty$.

\noindent
In the case b) of the lemma when $a_{k}=0$, we use the identity
$$
\widehat{G_{k}}(p)=\widehat{G_{k}}(0)+\int_{0}^{p}\frac{d\widehat{G_{k}}(s)}{ds}ds.
$$
Hence
\begin{equation}
\label{Gpib}  
\frac{\widehat{G_{k}}(p)}{p^{2}-ib_{k}p}=\frac{\widehat{G_{k}}(0)}{p(p-ib_{k})}+
\frac{\int_{0}^{p}\frac{d\widehat{G_{k}}(s)}{ds}ds}{p(p-ib_{k})}.
\end{equation}
By virtue of definition (\ref{ft}) of the standard Fourier transform, we
easily obtain the upper bound
\begin{equation}
\label{dgkdp}  
\Bigg|\frac{d\widehat{G_{k}}(p)}{dp}\Bigg|\leq \frac{1}{\sqrt{2\pi}}
\|xG_{k}(x)\|_{L^{1}({\mathbb R})}.
\end{equation}
We easily arrive at
$$
\Bigg|\frac{\int_{0}^{p}\frac{d\widehat{G_{k}}(s)}{ds}ds}{p(p-ib_{k})}\Bigg|\leq
\frac{\|xG_{k}(x)\|_{L^{1}({\mathbb R})}}{\sqrt{2\pi}|b_{k}|}<\infty
$$
via our assumptions. Thus, the expression in the left side
of (\ref{Gpib}) is bounded if and only if $\widehat{G_{k}}(0)$ vanishes, which
is equivalent to orthogonality relation (\ref{or1}).

\noindent
In the cases c) and d) of the lemma, we can write
$$
\frac{p^{2}\widehat{G_{k}}(p)}{p^{2}-a_{k}}=\widehat{G_{k}}(p)+
a_{k}\frac{\widehat{G_{k}}(p)}{p^{2}-a_{k}},
$$
so that
$\displaystyle{\frac{\widehat{G_{k}}(p)}{p^{2}-a_{k}}\in L^{\infty}({\mathbb R})}$
implies that
$\displaystyle{\frac{p^{2}\widehat{G_{k}}(p)}{p^{2}-a_{k}}\in
L^{\infty}({\mathbb R})}$ as well.

\noindent
To demonstrate the validity of the result of the part c) of our lemma, we
express
\begin{equation}
\label{gkpdelpm}
\frac{\widehat{G_{k}}(p)}{p^{2}-a_{k}}=\frac{\widehat{G_{k}}(p)}{p^{2}-a_{k}}
\chi_{I_{\delta_{k}}^{+}}+\frac{\widehat{G_{k}}(p)}{p^{2}-a_{k}}\chi_{I_{\delta_{k}}^{-}}+
\frac{\widehat{G_{k}}(p)}{p^{2}-a_{k}}\chi_{I_{\delta_{k}}^{c}}.
\end{equation}  
Here and below $\chi_{A}$ will denote the characteristic function of a set
$A\subseteq {\mathbb R}$ and $A^{c}$ will stand for its complement. The sets
$$
I_{\delta_{k}}^{+}:=\{p\in {\mathbb R} \ | \ \sqrt{a_{k}}-\delta_{k}<p<\sqrt{a_{K}}+
\delta_{k} \}, \quad
I_{\delta_{k}}^{-}:=\{p\in {\mathbb R} \ | \ -\sqrt{a_{k}}-\delta_{k}<p<-\sqrt{a_{k}}
+\delta_{k} \}
$$
with $0<\delta_{k}<\sqrt{a_{k}}$ and
$I_{\delta_{k}}:=I_{\delta_{k}}^{+}\cup I_{\delta_{k}}^{-}$.

\noindent
The third term in the right side of (\ref{gkpdelpm}) can be trivially bounded
from above in the absolute value by means of (\ref{inf1}) by
$\displaystyle{{1\over \sqrt{2\pi}\delta_{k}^{2}}\|G_{k}\|_{L^{1}({\mathbb R})}
<\infty}$.

\noindent
Clearly, we can write
$$
\widehat{G_{k}}(p)=\widehat{G_{k}}(\sqrt{a_{k}})+\int_{\sqrt{a_{k}}}^{p}\frac
{d\widehat{G_{k}}(s)}{ds}ds, \quad
\widehat{G_{k}}(p)=\widehat{G_{k}}(-\sqrt{a_{k}})+\int_{-\sqrt{a_{k}}}^{p}\frac
{d\widehat{G_{k}}(s)}{ds}ds.        
$$
Thus, the sum of the first two terms in the right side of (\ref{gkpdelpm}) is
given by
$$
\frac{\widehat{G_{k}}(\sqrt{a_{k}})}{p^{2}-a_{k}}\chi_{I_{\delta_{k}}^{+}}+
\frac{\widehat{G_{k}}(-\sqrt{a_{k}})}{p^{2}-a_{k}}\chi_{I_{\delta_{k}}^{-}}+
\frac{\int_{\sqrt{a_{k}}}^{p}\frac{d\widehat{G_{k}}(s)}{ds}ds}{p^{2}-a_{k}}
\chi_{I_{\delta_{k}}^{+}}+
\frac{\int_{-\sqrt{a_{k}}}^{p}\frac{d\widehat{G_{k}}(s)}{ds}ds}{p^{2}-a_{k}}
\chi_{I_{\delta_{k}}^{-}}.
$$
By virtue of (\ref{dgkdp}), we derive
$$
\Bigg|\frac{\int_{\sqrt{a_{k}}}^{p}\frac{d\widehat{G_{k}}(s)}{ds}ds}{p^{2}-a_{k}}
\chi_{I_{\delta_{k}}^{+}}\Bigg|\leq \frac{\|xG_{k}(x)\|_{L^{1}({\mathbb R})}}
{\sqrt{2\pi}(2\sqrt{a_{k}}-\delta_{k})}<\infty,
$$
$$
\Bigg|\frac{\int_{-\sqrt{a_{k}}}^{p}\frac{d\widehat{G_{k}}(s)}{ds}ds}{p^{2}-a_{k}}
\chi_{I_{\delta_{k}}^{-}}\Bigg|\leq \frac{\|xG_{k}(x)\|_{L^{1}({\mathbb R})}}
{\sqrt{2\pi}(2\sqrt{a_{k}}-\delta_{k})}<\infty.   
$$
Therefore, it remains to consider
\begin{equation}
\label{gkrakpm}  
\frac{\widehat{G_{k}}(\sqrt{a_{k}})}{p^{2}-a_{k}}\chi_{I_{\delta_{k}}^{+}}+
\frac{\widehat{G_{k}}(-\sqrt{a_{k}})}{p^{2}-a_{k}}\chi_{I_{\delta_{k}}^{-}}.
\end{equation}
Evidently, (\ref{gkrakpm}) is bounded if and only if
$\widehat{G_{k}}(\pm\sqrt{a_{k}})=0$. This is equivalent to the orthogonality
conditions (\ref{or12}).

\noindent
Finally, we turn our attention to the case d)  of the lemma when $a_{k}=0$, so
that
\begin{equation}
\label{gkp2}  
\frac{\widehat{G_{k}}(p)}{p^{2}}=\frac{\widehat{G_{k}}(p)}{p^{2}}
\chi_{\{|p|\leq 1\}}+\frac{\widehat{G_{k}}(p)}{p^{2}}
\chi_{\{|p|>1\}}.
\end{equation}
The second term in the right side of (\ref{gkp2}) can be easily estimated from
above in the absolute value as
$$
\Bigg|\frac{\widehat{G_{k}}(p)}{p^{2}}\chi_{\{|p|>1\}}\Bigg|\leq
\|\widehat{G_{k}}(p)\|_{L^{\infty}({\mathbb R})}<\infty
$$
via (\ref{inf1}). Obviously,
$$
\widehat{G_{k}}(p)=\widehat{G_{k}}(0)+p\frac{d\widehat{G_{k}}}{dp}(0)+
\int_{0}^{p}\Bigg(\int_{0}^{s}\frac{d^{2}\widehat{G_{k}}(q)}{dq^{2}}dq\Bigg)ds,
$$
such that the first term in right side of (\ref{gkp2}) equals to
\begin{equation}
\label{gkdgkp2}
\Bigg[  
\frac{\widehat{G_{k}}(0)}{p^{2}}+\frac{\frac{d\widehat{G_{k}}}{dp}(0)}{p}+
\frac{\int_{0}^{p}\Bigg(\int_{0}^{s}\frac{d^{2}\widehat{G_{k}}(q)}{dq^{2}}dq\Bigg)ds}
{p^{2}}\Bigg]\chi_{\{|p|\leq 1\}}.
\end{equation}
Using definition (\ref{ft}) of the standard Fourier transform,
we derive
$$
\Bigg|\frac{d^{2}\widehat{G_{k}}(p)}{dp^{2}}\Bigg|\leq \frac{1}{\sqrt{2\pi}}
\|x^{2}G_{k}(x)\|_{L^{1}({\mathbb R})},
$$
so that
$$
\Bigg|\frac{\int_{0}^{p}\Bigg(\int_{0}^{s}\frac{d^{2}\widehat{G_{k}}(q)}
{dq^{2}}dq\Bigg)ds}{p^{2}}\chi_{\{|p|\leq 1\}}\Bigg|\leq
\frac{\|x^{2}G_{k}(x)\|_{L^{1}({\mathbb R})}}{2\sqrt{2\pi}}<\infty
$$
as assumed. Hence, it remains to analyze
\begin{equation}
\label{gkdgk0}
\Bigg[  
\frac{\widehat{G_{k}}(0)}{p^{2}}+\frac{\frac{d\widehat{G_{k}}}{dp}(0)}{p}\Bigg]
\chi_{\{|p|\leq 1\}}.
\end{equation}
From definition (\ref{ft}) of the standard Fourier transform we deduce that
$$
\widehat{G_{k}}(0)=\frac{1}{\sqrt{2\pi}}(G_{k}(x), 1)_{L^{2}({\mathbb R})}  , \quad
\frac{d\widehat{G_{k}}}{dp}(0)=-\frac{i}{\sqrt{2\pi}}
(G_{k}(x), x)_{L^{2}({\mathbb R})},
$$
such that (\ref{gkdgk0}) is equal to
\begin{equation}
\label{gkdgk1}  
\frac{1}{\sqrt{2\pi}}\Bigg[\frac{(G_{k}(x), 1)_{L^{2}({\mathbb R})}}{p^{2}}-i
\frac{(G_{k}(x), x)_{L^{2}({\mathbb R})}}{p}\Bigg]\chi_{\{|p|\leq 1\}}.  
\end{equation}
Obviously, (\ref{gkdgk1}) belongs to $L^{\infty}({\mathbb R})$ if and only if
orthogonality relations (\ref{or13}) are valid. \hfill\lanbox

\bigskip

\noindent
Let the continuous function
$G_{k}(x): I\to {\mathbb R}, \ G_{k}(0)=G_{k}(2\pi)$. Its
Fourier transform on the finite interval is given by
\begin{equation}
\label{fti}
G_{k, n}:=\int_{0}^{2\pi}G_{k}(x)\frac{e^{-inx}}{\sqrt{2\pi}}dx, \quad n\in
{\mathbb Z},  
\end{equation}
so that
$\displaystyle{G_{k}(x)=\sum_{n=-\infty}^{\infty}G_{k, n}\frac{e^{inx}}{\sqrt{2\pi}}}$.
Obviously, the inequalities
\begin{equation}
\label{fubi}
\|G_{k, n}\|_{l^{\infty}}\leq \frac{1}{\sqrt{2\pi}}\|G_{k}(x)\|_{L^{1}(I)}, \quad
\|G_{k}(x)\|_{L^{1}(I)}\leq 2\pi\|G_{k}(x)\|_{C(I)}
\end{equation}
hold. 

\noindent
Similarly to the whole real line case, we define
\begin{equation}
\label{Nac}
{\cal N}_{a, \ b, \ k}:=\hbox{max}\Big\{
\Big\|\frac{G_{k, n}}{n^{2}-a_{k}-ib_{k}n}\Big\|_{l^{\infty}},
\quad
\Big\|\frac{n^{2}G_{k, n}}{n^{2}-a_{k}-ib_{k}n}\Big\|_{l^{\infty}}
\Big\}, \quad 1\leq k\leq K,
\end{equation}
\begin{equation}
\label{Mac}
{\cal M}_{a, \ k}:=\hbox{max}\Big\{
\Big\|\frac{G_{k, n}}{n^{2}-a_{k}}\Big\|_{l^{\infty}},
\quad
\Big\|\frac{n^{2}G_{k, n}}{n^{2}-a_{k}}\Big\|_{l^{\infty}}
\Big\}, \quad K+1\leq k\leq N
\end{equation}
with the constants $a_{k}\geq 0, \ b_{k}\in {\mathbb R}, \ b_{k}\neq 0$ and
$K\geq 2, \ N\geq 5$.
Let ${\cal N}_{0, \ b, \ k}$ stand for (\ref{Nac}) when $a_{k}=0$ and
${\cal M}_{0, \ k}$ denote (\ref{Mac}) when $a_{k}$ is trivial.
Under the conditions of Lemma A2 below, expressions (\ref{Nac}), (\ref{Mac})
will be finite.
This will enable us to introduce
\begin{equation}
\label{Nab2}
{\cal N}_{a, \ b}:=\hbox{max}_{1\leq k\leq K}{\cal N}_{a, \ b, \ k}<\infty,
\end{equation}
\begin{equation}
\label{Ma2}
{\cal M}_{a}:=\hbox{max}_{K+1\leq k\leq N}{\cal M}_{a, \ k}<\infty
\end{equation}
and
\begin{equation}
\label{Q2}
{\cal Q}:=\hbox{max}\{{\cal N}_{a, \ b}, \ {\cal M}_{a}\}.
\end{equation}  
Our final technical proposition is as follows.

\bigskip

\noindent
{\bf Lemma A2.} {\it Let $N\geq 5, \ K\geq 2, \
1\leq k\leq N, \ 1\leq l\leq K-1, \  K+1\leq r\leq q-1, \ r+1\leq q\leq N-1$,
the integral kernels
$G_{k}(x): I\to {\mathbb R}, \ G_{k}(x)\in C(I)$ and
$G_{k}(0)=G_{k}(2\pi)$ for all $1\leq k\leq N$.

\medskip

\noindent
a) Let $a_{k}>0, \ b_{k}\in {\mathbb R}, \ b_{k}\neq 0$ for $1\leq k\leq l$.
Then ${\cal N}_{a, \ b, \ k}<\infty$.

\medskip

\noindent
b) Let $a_{k}=0, \ b_{k}\in {\mathbb R}, \ b_{k}\neq 0$ for $l+1\leq k\leq K$.
Then ${\cal N}_{0, \ b, \ k}<\infty$ if and only if 
\begin{equation}
\label{or2}
(G_{k}(x),1)_{L^{2}(I)}=0.
\end{equation}

\medskip

\noindent
c) Let $a_{k}>0$ and $a_{k}\neq n^{2}, \ n\in {\mathbb Z}$ for
$K+1\leq k\leq r$. Then ${\cal M}_{a, \ k}<\infty$.

\medskip

\noindent
d) Let $a_{k}=n_{k}^{2}, \ n_{k}\in {\mathbb N}$ for $r+1\leq k\leq q$.
Then ${\cal M}_{a, \ k}<\infty$ if and only if 
\begin{equation}
\label{or21}
\Big(G_{k}(x), \frac{e^{\pm i n_{k}x}}{\sqrt{2\pi}}\Big)_{L^{2}(I)}=0  
\end{equation}  
holds.

\medskip

\noindent
e) Let $a_{k}=0$ for $q+1\leq k\leq N$. Then ${\cal M}_{0, \ k}<\infty$ if and
only if orthogonality condition (\ref{or2})
is valid.}

\medskip

\noindent
{\it Proof.} In both cases a) and b) of the lemma
$\displaystyle{\frac{G_{k, n}}{n^{2}-a_{k}-ib_{k}n}\in l^{\infty}}$ yields
$\displaystyle{\frac{n^{2}G_{k, n}}{n^{2}-a_{k}-ib_{k}n}\in l^{\infty}}$. Indeed,
$\displaystyle{\frac{n^{2}G_{k, n}}{n^{2}-a_{k}-ib_{k}n}}$ can be written as
\begin{equation}    
\label{Gkn3}
G_{k, n}+a_{k}\frac{G_{k, n}}{n^{2}-a_{k}-ib_{k}n}+ib_{k}
\frac{nG_{k, n}}{n^{2}-a_{k}-ib_{k}n}.  
\end{equation}
The first term in (\ref{Gkn3}) is bounded by means of (\ref{fubi}) under the
given conditions.
The third term in (\ref{Gkn3}) can be easily
estimated from above in the absolute value using (\ref{fubi}) as well as
$$
|b_{k}|\frac{|n||G_{k, n}|}
{\sqrt{(n^{2}-a_{k})^{2}+b_{k}^{2}n^{2}}}
\leq |G_{k, n}|\leq \sqrt{2\pi}\|G_{k}(x)\|_{C(I)}<\infty.
$$
Thus, $\displaystyle{\frac{G_{k, n}}{n^{2}-a_{k}-ib_{k}n}\in l^{\infty}}$ implies
the boundedness of
$\displaystyle{\frac{n^{2}G_{k, n}}{n^{2}-a_{k}-ib_{k}n}}$.

\noindent
To treat the case a) of the lemma, we need to consider
\begin{equation}
\label{Gknfr}
\frac{|G_{k, n}|}{\sqrt{(n^{2}-a_{k})^{2}+b_{k}^{2}n^{2}}}. 
\end{equation}
Obviously, the denominator in fraction (\ref{Gknfr}) can be estimated
from below by a positive constant. Let us apply (\ref{fubi})
to the numerator in (\ref{Gknfr}). Hence ${\cal N}_{a, \ b, \ k}<\infty$ in 
the situation a) of our lemma when $a_{k}>0$.

\noindent
Let us establish the validity of the result of the part b)
when $a_{k}$ vanishes. Evidently,
$\displaystyle{\frac{G_{k, n}}{n^{2}-ib_{k}n}\in l^{\infty}}$ if and only if
$G_{k, 0}=0$. This is equivalent to orthogonality relation (\ref{or2}).
In this case we easily obtain the upper bound
$$
\Bigg|\frac{G_{k, n}}{n^{2}-ib_{k}n}\Bigg|=\frac{|G_{k, n}|}
{|n|\sqrt{n^{2}+b_{k}^{2}}}\leq \sqrt{2\pi}\frac{\|G_{k}(x)\|_{C(I)}}{|b_{k}|}<\infty
$$
using (\ref{fubi}) along with our assumptions.

\noindent
In the cases c), d) and e) of the lemma we can express
\begin{equation}
\label{n2gkn}  
\frac{n^{2}G_{k, n}}{n^{2}-a_{k}}=G_{k, n}+a_{k}\frac{G_{k, n}}{n^{2}-a_{k}}.
\end{equation}
By virtue of (\ref{fubi}), the boundedness of  
$\displaystyle{\frac{G_{k, n}}{n^{2}-a_{k}}}$ yields
$\displaystyle{\frac{n^{2}G_{k, n}}{n^{2}-a_{k}}\in l^{\infty}}$.

\noindent
In the situation c) we have  
$\displaystyle{\frac{G_{k, n}}{n^{2}-a_{k}}\in l^{\infty}}$ since such expression
does not contain any singularities and the result of our lemma is obvious.

\noindent
In the part d) the quantity 
$\displaystyle{\frac{G_{k, n}}{n^{2}-n_{k}^{2}}\in l^{\infty}}$ if and only if
$G_{k, \pm n_{k}}=0$. This is equivalent to orthogonality conditions (\ref{or21}).

\noindent
In the case e) the expression 
$\displaystyle{\frac{G_{k, n}}{n^{2}}\in l^{\infty}}$ if and only if
$G_{k, 0}=0$, which is equivalent to orthogonality relation (\ref{or2}).
\hfill\lanbox 

\bigskip


\bigskip

\section*{Acknowledgement} V. V. is grateful to Israel Michael Sigal
for the partial support by the NSERC grant NA 7901.

\bigskip


\end{document}